\documentclass[11pt,a4paper]{amsart} 

\usepackage{latexsym}

\usepackage[a4paper , lmargin = {2.5cm} , rmargin = {2.5cm} , tmargin = {2.5cm} , bmargin = {2.5cm} ]{geometry} 

\usepackage{geometry}
\usepackage{pinlabel}         
\usepackage{graphicx}
\usepackage{amssymb, amsthm}
\usepackage{epstopdf}
\usepackage{romannum}
\usepackage{tikz}
\usepackage{subfig}
\usepackage{MnSymbol}

  	\newcommand{\Z}{\ensuremath{\mathbb{Z}}}   
  	   
  	\newcommand{\R}{\ensuremath{\mathbb{R}}} 

   	\def\Aut{{\rm{Aut}}$(F_n)$}
        \def\Autn{{\rm{Aut}}(F_n)}
        \def\Autm{{\rm{Aut}}(F_{n-1})}
   	\def\SAut{{\rm{SAut}}$(F_n)$}
        \def\SAutn{{\rm{SAut}}(F_n)}
	\def\GL{{\rm{GL}}}
	\def\SL{{\rm{SL}}}
        \def\Isom{{\rm{Isom}}$(X)$}
	\def\isom{{\rm{Isom}}(X)}
	\def\ord{{\rm{ord}}}
	\def\im{{\rm{im}}}
	\def\H{{\rm{H}}}
	\def\Fix{{\rm{Fix}}}
 	
	\def\FAd{{\rm{FA}}$_{d}$}
	
	\def\CAT{{\rm{CAT}}$(0)$}
	\def\Sym{{\rm{Sym}}(n)}

	\def\F{{\rm{F}}}

\theoremstyle{plain}
\newtheorem*{NewTheorem}{Theorem}	
\newtheorem*{NewTheoremA}{Theorem A}
\newtheorem*{NewTheoremB}{Theorem B}
\newtheorem*{NewCorollaryC}{Corollary C}
\newtheorem*{NewCorollaryD}{Corollary D}
\newtheorem*{Remark}{Remark}
\newtheorem*{FarbF}{Farb's Fixed Point Criterion}

\newtheorem{theorem}{Theorem}[section]

\newtheorem{corollary}[theorem]{Corollary}
\newtheorem{proposition}[theorem]{Proposition}
\newtheorem{definition}[theorem]{Definition}

\newtheorem{convention}[theorem]{Convention}
\newtheorem{HellyRTheorem}[theorem]{Helly's classical Theorem}
\newtheorem{HellyOpenTheorem}[theorem]{Helly's Theorem for open convex subspaces of a $\mathbf{CAT(0)}$ space}
\newtheorem{HellyClosedTheorem}[theorem]{Helly's Theorem for closed convex subspaces of a $\mathbf{CAT(0)}$ space}
\newtheorem{Farb}[theorem]{Farb's Fixed Point Criterion}
\newtheorem{remark}[theorem]{Remark}

\newenvironment{pf}{\par\medskip\noindent\textit{Proof.}~}{\hfill $\square$\par\medskip}

\title[Fixed points for actions of Aut$(F_{n})$ on \CAT\ spaces]{Fixed points for actions of Aut\boldmath$(F_{n})$\unboldmath\ on CAT\boldmath$(0)$\unboldmath\ spaces}
\author{Olga Varghese}
\date{\today}
\address{Olga Varghese\\
Department of Mathematics\\
M\"unster University\\ 
Einsteinstra\ss e 62\\
48149 M\"unster (Germany)}
\email{olga.varghese@uni-muenster.de}

\begin{document}
\pagenumbering{arabic}
\begin{abstract}
For $n\geq 4$ we discuss questions concerning global fixed points for isometric actions of \Aut, the automorphism group of a free group of rank $n$, on complete \CAT\ spaces. We prove that whenever \Aut\ acts by isometries on complete $d$-dimensional \CAT\ space with $d<2\left\lfloor\frac{n}{4}\right\rfloor-1$, then it must fix a point. This property has implications for irreducible representations of \Aut, which are also presented here. For \SAut, the unique subgroup of index two in \Aut, we obtain similar results.
\end{abstract}

\maketitle

\section{Introduction}
In the mathematical world, this article is located in the area of geometric group theory, a field at the intersection of algebra, geometry and topology. Geometric group theory studies the interaction between algebraic and geometric properties of groups. One~is interested to understand on which 'nice' geometric spaces a given  group can act in a reasonable way and how geometric properties of these spaces are reflected in the algebraic structure of the group.
Here, the spaces will  be \CAT\ metric spaces,  while the groups will be \Aut\ and \SAut. 
The questions we shall investigate are concerned with fixed point properties and the representation theory of these groups.

More precisely, let $\mathbb{Z}^n$ be the free abelian group and $F_{n}$ the free group of rank $n$. One goal for a group theorist is to unterstand the structure of their automorphism groups, 
${\rm GL}_n(\mathbb{Z})$ resp. \Aut.
The abelianization map $F_{n}\twoheadrightarrow \Z^{n}$ gives a natural epimorphism 
$\Autn\twoheadrightarrow\GL_{n}(\Z).$
The special automorphism group of $F_n$, which we will denote by \SAut, is defined as the preimage of ${\rm SL}_{n}(\Z)$ under this map. 
Much of the work on \Aut\ and \SAut\ is motivated by the idea that ${\rm GL}_n(\mathbb{Z})$ and \Aut\ resp. ${\rm SL}_n(\Z)$ and \SAut\  should have many properties in common. 
Here we follow this idea and present analogies between these groups with respect to fixed point properties.

Let $\mathcal{X}$ be a class of metric spaces. A group $G$ is said to have property $\F\mathcal{X}$ if any action of $G$ by isometries on any member of $\mathcal{X}$ has a fixed point. Let $\mathcal{A}$ be the class of simplicial trees, $\mathcal{A}_d$ the class of complete \CAT\ spaces of covering dimension~$d$ and $\mathcal{A}_{*}$ the class of finite dimensional complete \CAT\ spaces. 

The starting point for our investigation is the study of group actions on simplicial trees which was initiated by \textsc{Serre}, \cite{Serre1}, \cite{Serre2}. He proved that ${\rm GL}_n(\Z)$ and~${\rm SL}_n(\Z)$ have property $\F\mathcal{A}$ for $n\geq 3$. Regarding \Aut\ and \SAut, \textsc{Bogopolski} was the first to prove that these groups also have property $\F\mathcal{A}$, see \cite{BogopolskiFA}.

A slight generalization of the class of simplicial trees is given by the class of metric trees, which we will denote by $\mathcal{R}$. Different methods were developed by \textsc{Culler} and \textsc{Vogtmann} and later by \textsc{Bridson} to prove that \Aut\ and \SAut\ have property $\F\mathcal{R}$, \cite{BridsonFA1}, \cite{VogtmannFA}. We obtain the fixed point property of \Aut\ and \SAut\ for a much larger class of higher dimensional complete \CAT\ spaces.

We present two results Theorems A and B regarding property $\F\mathcal{A}_{*}$ for the groups \Aut\ and \SAut. Using \textsc{Bridson's} and \textsc{Farb's} techniques from \cite{BridsonFA1} and \cite{Farb}, we prove:
\begin{NewTheoremA}
If $n\geq 4$ and $d<{\rm min}\left\{k\left\lfloor\frac{n}{k+2}\right\rfloor\mid k=2,\ldots, d+1\right\}$, then \Aut\ has property $\F\mathcal{A}_{d}$. In particular, if $n\geq 4$ and $d<2\left\lfloor\frac{n}{4}\right\rfloor-1$, then \Aut\ has property~$\F\mathcal{A}_{d}$. 
\end{NewTheoremA}

\begin{NewTheoremB}
If $n\geq 5$ and $d<{\rm min}\left\{k\left\lfloor\frac{n-1}{k+2}\right\rfloor\mid k=2,\ldots,d+1\right\}$, then \SAut\ has~property $\F\mathcal{A}_{d}$.
In particular, if $n\geq 5$ and $d<2\left\lfloor\frac{n-1}{4}\right\rfloor-1$, then \SAut\ has property~$\F\mathcal{A}_{d}$.
\end{NewTheoremB}

Our proofs of Theorems A and B involve three ingredients. First we construct a generating set of \Aut\ such that each pair of its elements generates a finite subgroup. Next,
we need an extended version of Helly's Theorem for higher dimensional~\CAT\ spaces.
\begin{NewTheorem}
 Let $X$ be a $d$-dimensional complete \CAT\ space and $\mathcal{S}$ a finite family of non-empty closed convex subspaces. If the intersection of each $(d+1)$-elements of $\mathcal{S}$ is non-empty, then~$\bigcap\mathcal{S}$ is non-empty.
\end{NewTheorem}
There exist several variations of this theorem in the literature, e.g. for finite families of convex open resp. closed subsets of a \CAT\ space, see \cite[3.2]{MappingClassFA}, \cite[2]{Debrunner},  \cite[3.2]{Farb} and \cite[5.3]{Kleiner}. Here we include a complete proof for the case of a finite family of closed convex subspaces.

Our main technique in the proofs of Theorems A and B is based on the following corollary.  Indeed, it was \textsc{Farb} who discovered the connection between Helly's Theorem and the combinatorics of generating sets for a large class of groups.
\begin{FarbF}
\label{HellyGroup}
Let $G$ be a group, $Y$ a finite generating set of~$G$ and~$X$ a complete $d$-dimensional \CAT\ space. If
$\Phi:G\rightarrow\isom$ is a homomorphism such that each $(d+1)$-element subset of $Y$ has a fixed point in $X$, then $G$ has a fixed point in $X$.
\end{FarbF}
\textsc{Farb} used this criterion in \cite{Farb} to obtain sharp results on property $\F\mathcal{A}_d$ for various groups. For example, he proved that  ${\rm SL}_n(\mathbb{Z}[1/p])$ has property $\F\mathcal{A}_{n-2}$ for semisimple actions, but not $\F\mathcal{A}_{n-1}$, since it acts without a global fixed point on the affine building for ${\rm SL}_n(\mathbb{Q}_p)$.

In a third step, we combine the extended version of Helly's Theorem with the following theorem by \textsc{Bridson}  to prove our results.

\begin{NewTheorem}(\cite[3.6]{MappingClassFA})
 Let $k$ and $l$ be in $\mathbb{N}_{>0}$ and let $X$ be a complete $d$-dimensional~\CAT\ space with $d<k\cdot l$. Let $S$ be a subset of \Isom\ and let $S_{1},\ldots, S_{l}$ be conjugates of $S$ such that $[S_{i}, S_{j}]=1$ for $i,j=1,\ldots, l,\ i\neq j$. If each $k$-element subset of $S$ has a fixed point in $X$, then each finite subset of $S$ has a fixed point in $X$.
\end{NewTheorem}

Property $\F\mathcal{A}_{d}$ strongly affects the representation theory of groups. The following result by \textsc{Farb}, partially based on work by \textsc{Bass}, illustrates this fact.

\begin{NewTheorem}(\cite[1.8]{Farb})
Let $K$ be an algebraically closed field and let $G$ be a group. If $G$ has property $\F\mathcal{A}_d$, then
there are only finitely many conjugacy classes of irreducible representations
\[
 \rho:G\rightarrow{\rm GL}_{d+1}(K).
\]
\end{NewTheorem}

As an application of our Theorems A and B, we obtain the following similar results for the representation theory of \Aut\ and \SAut. 
\begin{NewCorollaryC}
Let $K$ be an algebraically closed field. If $n\geq 4$ and $d\leq 2\left\lfloor\frac{n}{4}\right\rfloor-1$, then
there are only finitely many conjugacy classes of irreducible representations
\[
\rho:\Autn\rightarrow{\rm GL}_{d}(K).
\]
\end{NewCorollaryC}

\begin{NewCorollaryD}
Let $K$ be an algebraically closed field. If $n\geq 5$ and $d\leq 2\left\lfloor\frac{n-1}{4}\right\rfloor-1$, then
there are only finitely many conjugacy classes of irreducible representations
\[
\rho:\SAutn\rightarrow{\rm SL}_{d}(K).
\]
\end{NewCorollaryD}

\begin{Remark}
A better bound for the complex representations of \Aut\ is proved in \cite[3.1,3.2]{Rapinchuk}.
If $n\geq 3$ and $d\leq2\cdot n-2$, then there are only finitely many conjugacy classes of irreducible representations
\[
 \rho:\Autn\rightarrow\GL_{d}(\mathbb{C}).
\]
With \textsc{Bridson's} and \textsc{Vogtmann's} techniques from \cite[1.1]{VogtmannSpheres} one can prove that the linear representations of \SAut\ are very rigid.
Let $K$ be a field of characteristic not equal to two and let
\[
\rho:\SAutn\rightarrow\SL_{d}(K)
\]
be a homomorphism. If $n\geq 3$ and $d<n$, then $\rho$ is trivial.
In particular, if $n\geq 3$, then \Aut\ has only finitely many conjugacy classes of irreducible representations in any dimension $\leq n-1$.
\end{Remark}

\section{A generating set of \Aut}

The purpose of this section is to construct a generating set of the group \Aut\ such that each pair of its elements generates a finite subgroup. 

Although it seems awkward at first glance, it is convenient and standard to work with the right action of \Aut\ on $F_n$.
\begin{convention}
For $\alpha, \beta$ in $\Autn$ the automorphism $\alpha\beta$ is the composite
where $\alpha$ acts before $\beta$.
\end{convention}

Let us first introduce a notations for some elements of \Aut.   We define the {\bf right Nielsen automorphism} $\rho_{ij}$,  involutions $(x_{i},x_{j})$ and $e_{i}$ for $i, j=1,\ldots, n$,~$i\neq j$ as follows:
\[
\begin{matrix}
\rho_{ij}(x_{k}):=\begin{cases} x_{i}x_{j} & \mbox{if $k=i$,}  \\ 
				     x_{k} & \mbox{if $k\neq i$.}   
\end{cases}
&
(x_{i},x_{j})(x_{k}):=\begin{cases} x_{j} & \mbox{if $k=i$,}  \\ 
			     x_{i} & \mbox{if $k=j$,}  \\ 
			     x_{k} & \mbox{if $k\neq i,$ $j$.}  
\end{cases}
&
e_{i}(x_{k}):=\begin{cases} x^{-1}_{i} & \mbox{if $k=i$,}  \\
				 x_{k} & \mbox{if $k\neq i$}.   
\end{cases}
\end{matrix}
\]

It is easy to see that the image of $X=\left\{x_1,\ldots, x_n\right\}$ under any of these maps is another basis of $F_n$, therefore these elements  are automorphisms. It was proven by \textsc{Nielsen} in \cite[p.~173]{Nielsen}) that for~$n\geq 3$ the group \Aut\  is generated by the set 
\[
Y_{1}:=\left\{\rho_{12}, \ e_{1},\  (x_{1},x_{2}),\ (x_{1},x_{2},\ldots,x_{n})\right\},
\]
where $(x_{1},x_{2},\ldots,x_{n})$ is equal to $(x_{n-1}, x_{n})(x_{n-2}, x_{n-1})\ldots (x_{1},x_{2})$.

Our strategy in this section is to modify the set $Y_{1}$ such that each pair of elements in the new generating set generates a finite group, compare \cite[1.1, 1.2]{BridsonFA1}.

\begin{proposition}
\label{GenAut}
Let $n\geq 3$. 
\begin{enumerate}
 \item[$(i)$] The group \Aut\ is generated by
\[
Y_{2}:=\left\{(x_{1}, x_{2})e_{1}e_{2},\  (x_{2},x_{3})e_{1},\ (x_{i}, x_{i+1}),\ e_{2}\rho_{12}, \ e_{n}\mid i=3,\ldots, n-1\right\}.
\]
\item[$(ii)$] The subgroup generated by $Y_{2}-\left\{e_{2}\rho_{12}\right\}$ is finite.
\item[$(iii)$] For $\alpha, \beta$ in $Y_{2}$ the subgroup generated by $\left\{\alpha, \beta\right\}$ is finite.
\end{enumerate}
\end{proposition}
\begin{pf}
Let us denote by $\Sigma_n\subseteq\Autn$ the group of automorphisms which permute the basis $X$. The conjugation by $\sigma\in\Sigma_n$ sends $e_{i}$ to $e_{\sigma(i)}$:
$\sigma^{-1}e_i\sigma=e_{\sigma(i)}$,
therefore \Aut\ is generated by the set $\left\{\rho_{12},e_{n},\Sigma_n\right\}$.
It is a well-known result that the group $\Sigma_n$ is generated by the involutions $(x_{i},x_{i+1})$ with $i=1,\ldots, n-1$. We can further replace $\rho_{12}$ by the involution $e_{2}\rho_{12}$ and we obtain the following generating set of \Aut: $\left\{ e_{2}\rho_{12}, e_{n}, (x_{i}, x_{i+1})\mid i=1,\ldots,n-1\right\}$. To see that $Y_2$ is a generating set of \Aut, we must show that the involutions $(x_{1}, x_{2})$ and $(x_{2}, x_{3})$ are in $\langle Y_{2}\rangle$. First we show this results for $n=3$. We have
$e_2=\underbrace{(x_2, x_3)e_1}_{\in Y_2}\underbrace{e_3}_{\in Y_2}\underbrace{(x_2, x_3)e_1}_{\in Y_2}\in\langle Y_2\rangle$
and therefore
$(x_1, x_3)=\underbrace{(x_2,x_3)e_1}_{\in Y_2}\underbrace{(x_1, x_2)e_1e_2}_{\in Y_2}\underbrace{e_2}_{\in\langle Y_2\rangle}\underbrace{(x_2,x_3)e_1}_{\in Y_2}\underbrace{e_3}_{\in Y_2}$ is contained in $\langle Y_2\rangle$.
Using $(x_1, x_3, x_2)=\underbrace{(x_1, x_2)e_1e_2}_{\in Y_2}\underbrace{e_2}_{\in\langle Y_2\rangle}\underbrace{(x_2, x_3)e_1}_{\in Y_2}\in\langle Y_2\rangle$
we obtain
\begin{align*}
(x_1, x_2)&=\underbrace{(x_1, x_3)}_{\in\langle Y_2\rangle}\underbrace{(x_1, x_3, x_2)}_{\in\langle Y_2\rangle}\in \langle Y_2\rangle,\\
(x_2, x_3)&=\underbrace{(x_1,x_2)}_{\in\langle Y_2\rangle}\underbrace{(x_1, x_3, x_2)}_{\in\langle Y_2\rangle}\in\langle Y_2\rangle.
\end{align*}
For $n\geq 4$ we have $e_3=\underbrace{(x_3, x_n)}_{\in\langle Y_2\rangle}\underbrace{e_n}_{\in Y_2}\underbrace{(x_3,x_n)}_{\in\langle Y_2\rangle}\in\langle Y_2\rangle$. The same arguments as above show that the involutions $(x_1, x_2)$ and $(x_2, x_3)$ are contained in $\langle Y_2\rangle$.
This finishes the proof of statement $(i)$.

It easy to verify that the subgroup $\langle Y_{2}-\left\{e_{2}\rho_{12}\right\} \rangle$ of \Aut\ is isomorphic to the semidirect product $\Sym\ltimes\mathbb{Z}^{n}_{2}$, where we denote by $\mathbb{Z}_2$ the cyclic group of order $2$, and therefore finite.

Now we prove the last statement of the proposition. If $\left\{\alpha, \beta\right\}$ is a subset of $Y_{2}-\left\{e_{2}\rho_{12}\right\}$ then the statement is obvious. Otherwise we compute the order of $e_{2}\rho_{12}\alpha$ for $\alpha\in Y_2$ in detail.
The involution $e_{2}\rho_{12}$ commutes with $(x_{i},x_{i+1})$ for $i=3,\ldots,n$ and with $e_{n}$. It follows that the order of $e_{2}\rho_{12}(x_{i},x_{i+1})$ and $e_{2}\rho_{12}e_{n}$ is equal to two and therefore the subgroups
$\langle\left\{e_{2}\rho_{12}, (x_{i},x_{i+1})\right\}\rangle$ for $i=3,\ldots,n$ and  $\langle\left\{e_2\rho_{12}, e_{n}\right\}\rangle$ are isomorphic to  
$\mathbb{Z}_2\times\mathbb{Z}_2$.
The order of $e_{2}\rho_{12}(x_{1},x_{2})e_{1}e_{2}$ is equal to $3$,

it follows that the dihedral group ${\rm D}_{3}$ of order $6$ has an epimorphism onto
 $\langle\left\{e_{2}\rho_{12},(x_{1},x_{2})e_{1}e_{2}\right\}\rangle$ and this group is therefore finite.
The order of $e_{2}\rho_{12}(x_{2},x_{3})e_{1}$ is equal to $4$, 
 and it follows that the dihedral group ${\rm D}_{4}$ of order $8$ has an epimorphism onto 
$\langle\left\{e_{2}\rho_{12},(x_{2},x_{3})e_{1}\right\}\rangle$ and this group is therefore finite.
\end{pf}

\section{A generating set of \SAut}

In this section we construct a generating set for the group \SAut\ with the same finiteness property as the set $Y_2$ in the previous section. 

For $i, j=1,\ldots,n$, $i\neq j$ we define the {\bf left Nielsen automorphism }$\lambda_{ij}$ as follows:
\[
\begin{matrix}
\lambda_{ij}(x_{k}):=\begin{cases} x_{j}x_{i} & \mbox{if $k=i$,}  \\ 
				        x_{k} & \mbox{if $k\neq i$}.  
\end{cases}
\end{matrix}
\]
It is known that  the group \SAut\ is generated by $\left\{\rho_{ij},\lambda_{ij}\mid i, j=1,\ldots,n, \ i\neq j\right\}$ for $n\geq 3$, see \cite[2.8]{Gersten}. An easy calculation shows that the commutator of $\rho_{ij}$ and~$\rho_{jk}$ is equal to $\rho_{ik}$ and the commutator of $\lambda_{ij}$ and $\lambda_{jk}$ is equal to $\lambda_{ik}$ for $i,j,k=1,\ldots,n$ distinct, therefore  \SAut\ is generated by the set 
\[
Y_{3}=\left\{\rho_{i(i+1)},\  \rho_{n1},\ \lambda_{i(i+1)},\ \lambda_{n1}\mid i=1,\ldots, n-1\right\}.
\]

Our strategy in this section is to modify the set $Y_3$ to obtain a new generating set of \SAut\ which will have the additional property that each group generated by any two of its elements is finite.
\begin{proposition}
\label{finite2}
Let $n\geq 4$. 
\begin{enumerate}
 \item[$(i)$] The group \SAut\ is generated by 
\[
Y_{4}:=\left\{(x_{1}, x_{2})e_{1}e_{2}e_{3}, \ (x_{2},x_{3})e_{1},\ (x_{i}, x_{i+1})e_{i}, \ e_{2}e_{4}\rho_{12}, \ e_{3}e_{4}\mid i=3,\ldots, n-1\right\}.
\]
\item[$(ii)$] The subgroup generated by $Y_{4}-\left\{e_{2}e_{4}\rho_{12}\right\}$ is finite.
\item[$(iii)$] For $\alpha, \beta$ in $Y_{4}$ the subgroup generated by $\left\{\alpha, \beta\right\}$ is finite.
\end{enumerate}
\end{proposition}
\begin{pf}
Using the relation $e_{i}e_{j}\rho_{ij}e_{j}e_{i}=\lambda_{ij}$ for $i, j=1,\ldots,n$ with $i\neq j$ we obtain
$\SAutn=\langle\left\{\rho_{i(i+1)},\rho_{n1}, e_{i}e_{i+1}, e_{n}e_{1}\mid i=1,\ldots,n-1\right\}\rangle$. 
As a next step in the proof, we claim that \SAut\ is generated by the set 
\[
Y'=\left\{(x_{1}, x_{2})e_{1}e_{2}e_{3},\ (x_{2},x_{3})e_{1},\ (x_{i}, x_{i+1})e_{3},\ e_{2}e_{4}\rho_{12}, \ e_{3}e_{4}\mid i=3,\ldots, n-1\right\}.
\]
The element $e_{2}e_{4}=\underbrace{(x_{2},x_{3})e_{1}}_{\in Y'}\underbrace{e_{3}e_{4}}_{\in Y'}\underbrace{e_{1}(x_{2},x_{3})}_{\in \langle Y'\rangle}$ is contained in $\langle Y'\rangle$, therefore 
\[
\rho_{12}=\underbrace{e_{4}e_{2}}_{\in\langle Y'\rangle}\underbrace{e_{2}e_{4}\rho_{12}}_{\in Y'}\in\langle Y'\rangle.
\] 
From the relation $e_{2}e_{3}=\underbrace{e_{2}e_{4}}_{\in\langle Y'\rangle}\underbrace{e_{3}e_{4}}_{\in\langle Y'\rangle}\in\langle Y'\rangle$ we see that
\[
(x_{2},x_{3})(x_{1},x_{2})=\underbrace{(x_{2},x_{3})e_{1}}_{\in Y'}\underbrace{e_{2}e_{3}}_{\in\langle Y'\rangle}\underbrace{e_{3}e_{2}e_{1}(x_{1},x_{2})}_{\in\langle Y'\rangle}\in\langle Y'\rangle.
\]
Using 
$e_{1}e_{4}=\underbrace{(x_{1},x_{2})e_{1}e_{2}e_{3}}_{\in Y'}\underbrace{e_{2}e_{4}}_{\in\langle Y'\rangle}\underbrace{e_{3}e_{2}e_{1}(x_{1},x_{2})}_{\in\langle Y'\rangle}\in\langle Y'\rangle$ and
$e_{3}e_{1}=\underbrace{e_{1}e_{4}}_{\in\langle Y'\rangle}\underbrace{e_{3}e_4}_{\in\langle Y'\rangle}\in\langle Y'\rangle$ we see that 
\[
(x_{3},x_{4})(x_{2},x_{3})=\underbrace{(x_{3},x_{4})e_{3}}_{\in Y'}\underbrace{e_{3}e_{1}}_{\in\langle Y'\rangle}\underbrace{e_{1}(x_{2},x_{3})}_{\in\langle Y'\rangle}\in\langle Y'\rangle.
\]
Now we show that the element $e_{3}(x_{1}, x_{n})$ is contained in $\langle Y'\rangle$. We consider the relation  $e_{1}e_{2}=\underbrace{(x_{2},x_{3})e_{1}}_{\in Y'}\underbrace{e_{1}e_{3}}_{\in\langle Y'\rangle}\underbrace{e_{1}(x_{2},x_{3})}_{\in\langle Y'\rangle}\in\langle Y'\rangle$, therefore $(x_1, x_2)e_3=\underbrace{(x_1, x_2)e_1e_2e_3}_{\in Y'}\underbrace{e_1e_2}_{\in\langle Y'\rangle}$ is contained in $\langle Y'\rangle$.
If $n$ is odd, then we have 
\begin{multline*}
e_{3}(x_{1},x_{n})=\underbrace{e_{3}(x_{1},x_{2})}_{\in\langle Y'\rangle}\underbrace{(x_{2},x_{3})e_{3}}_{\in\langle Y'\rangle}\underbrace{e_{3}(x_{3},x_{4})}_{\in\langle Y'\rangle}\\
\ldots\underbrace{e_{3}(x_{n-2},x_{n-1})}_{\in\langle Y'\rangle}\underbrace{(x_{n-1},x_{n})e_{3}}_{\in Y'}\underbrace{e_{3}(x_{n-2},x_{n-1})}_{\in\langle Y'\rangle}
\\
\ldots\underbrace{(x_{2},x_{3})e_{3}}_{\in\langle Y'\rangle}\underbrace{e_{3}(x_{1},x_{2})}_{\in\langle Y'\rangle}
\end{multline*}
and if $n$ is even, then
\begin{multline*}
(x_{1},x_{n})e_{3}=\underbrace{(x_{1},x_{2})e_{3}}_{\in\langle Y'\rangle}\underbrace{e_{3}(x_{2},x_{3})}_{\in\langle Y'\rangle}\underbrace{(x_{3},x_{4})e_{3}}_{\in Y'}\\
\ldots\underbrace{e_{3}(x_{n-2},x_{n-1})}_{\in\langle Y'\rangle}\underbrace{(x_{n-1},x_{n})e_{3}}_{\in Y'}\underbrace{e_{3}(x_{n-2},x_{n-1})}_{\in\langle Y'\rangle}\\
\ldots\underbrace{e_{3}(x_{2},x_{3})}_{\in\langle Y'\rangle}\underbrace{(x_{1},x_{2})e_{3}}_{\in\langle Y'\rangle}. 
\end{multline*}
Using $e_{3}(x_{1}, x_{n})\in\langle Y'\rangle$ and $(x_{n-1},x_{n})e_{3}\in Y'$ we obtain $(x_{1}, x_{n})(x_{n-1}, x_{n})\in\langle Y'\rangle$.
From the relations
\begin{align*}
 \rho_{(i+1)(i+2)}&=\underbrace{(x_{i},x_{i+1})(x_{i+1},x_{i+2})}_{\in\langle Y'\rangle}\underbrace{\rho_{i,i+1}}_{\in\langle Y'\rangle}\underbrace{(x_{i+1},x_{i+2})(x_{i},x_{i+1})}_{\in\langle Y'\rangle},\\
    e_{i+1}e_{i+2}&=\underbrace{(x_{i},x_{i+1})(x_{i+1},x_{i+2})}_{\in\langle Y'\rangle}\underbrace{e_{i}e_{i+1}}_{\in\langle Y'\rangle}\underbrace{(x_{i+1},x_{i+2})(x_{i},x_{i+1})}_{\in\langle Y'\rangle}\\
\end{align*}
we see that $\rho_{i(i+1)}, e_{i}e_{i+1}$ are in $\langle Y'\rangle$ for $i=1,\ldots,n-1$ and using relations
\begin{align*}
  \rho_{n1}&=\underbrace{(x_{n-1},x_{n})(x_{1},x_{n})}_{\in\langle Y'\rangle}\underbrace{\rho_{(n-1)n}}_{\in\langle Y'\rangle}\underbrace{(x_{1},x_{n})(x_{n-1},x_{n})}_{\in\langle Y'\rangle},\\
 e_{n}e_{1}&=\underbrace{(x_{n-1},x_{n})(x_{1},x_{n})}_{\in\langle Y'\rangle}\underbrace{e_{n-1}e_{n}}_{\in\langle Y'\rangle}\underbrace{(x_{1},x_{n})(x_{n-1},x_{n})}_{\in\langle Y'\rangle}\\
\end{align*}
we obtain $\rho_{n1},e_{n}e_{1}\in\langle Y'\rangle$ and therefore $Y'$ is a generating set of \SAut.
Now we show that $(x_{i},x_{i+1})e_{3}$ is contained in $\langle Y_{4}\rangle$ for $i=4,\ldots,n-1$.
We have the relations
\begin{align*}
e_{5}e_{3}&=\underbrace{(x_{4},x_{5})e_{4}}_{\in Y_{4}}\underbrace{e_{4}e_{3}}_{\in\langle Y_{4}\rangle}\underbrace{e_{4}(x_{4},x_{5})}_{\in\langle Y_{4}\rangle}\in\langle Y_{4}\rangle,\\
e_{6}e_{3}&=\underbrace{(x_{5},x_{6})e_{5}}_{\in Y_{4}}\underbrace{e_{5}e_{3}}_{\in\langle Y_{4}\rangle}\underbrace{e_{5}(x_{5},x_{6})}_{\in\langle Y_{4}\rangle}\in\langle Y_{4}\rangle,\\
\ldots\\
e_{n-1}e_{3}&=\underbrace{(x_{n-2},x_{n-1})e_{n-2}}_{\in Y_{4}}\underbrace{e_{n-2}e_{3}}_{\in\langle Y_{4}\rangle}\underbrace{e_{n-2}(x_{n-2},x_{n-1})}_{\in\langle Y_{4}\rangle}\in\langle Y_{4}\rangle
\end{align*}
and we see that $e_{i}e_{3}\in\langle Y_{4}\rangle$ for $i=4,\ldots,n-1$ and
$(x_{i},x_{i+1})e_{3}=\underbrace{(x_{i},x_{i+1})e_{i}}_{\in Y_{4}}\underbrace{e_{i}e_{3}}_{\in\langle Y_{4}\rangle}\in\langle Y_{4}\rangle$, hence $\SAutn=\langle Y_{4}\rangle$.
 
Now we prove the second statement of the proposition. It is easy to verify that the subgroup $\langle Y'-\left\{e_{2}e_{4}\rho_{12}\right\}\rangle$ of \SAut\ is isomorphic to a subgroup of the semidirect product $\Sym\ltimes\mathbb{Z}^{n}_{2}$ and therefore finite.

For the proof of the last statement of the proposition we note that the elements in~$Y_{4}$ have finite order. If $\left\{\alpha, \beta\right\}$ is a subset of $Y_{4}-\left\{e_{2}e_{4}\rho_{12}\right\}$, then the statement is obvious. Otherwise we consider the subset $\left\{e_{2}e_{4}\rho_{12}, \alpha\right\}$ for $\alpha\in Y_4$. If the commutator of $e_{2}e_{4}\rho_{12}$ and $\alpha$ is equal to one, then the subgroup 
$\langle\left\{e_{2}e_{4}\rho_{12}, \alpha\right\}\rangle$ is finite. If $e_{2}e_{4}\rho_{12}$ does not commute with $\alpha$, then 
$\alpha\in\left\{(x_{1},x_{2})e_{1}e_{2}e_{3},(x_{2},x_{3})e_{1},(x_{3},x_{4})e_{3}, (x_{4},x_{5})e_{4}\right\}$. We note that
$\ord((x_{1},x_{2})e_{1}e_{2}e_{3})=\ord((x_{2},x_{3})e_{1})=2$, 
$\ord(e_{2}e_{4}\rho_{12}(x_{1},x_{2})e_{1}e_{2}e_{3})=6$ and $\ord(e_{2}e_{4}\rho_{12}(x_{2},x_{3})e_{1})=4$. It follows that the subgroups
$\langle\left\{e_{2}e_{4}\rho_{12}, (x_{1},x_{2})e_{1}e_{2}e_{3}\right\}\rangle$ and $\langle\left\{e_{2}e_{4}\rho_{12}, (x_{2},x_{3})e_{1}\right\}\rangle$ are finite. If $\alpha$ is equal to $(x_{3},x_{4})e_{3}$ or $(x_{4},x_{5})e_{4}$, then the dihedral group ${\rm D}_{4}:=\langle r,s\mid r^{4}=s^{2}=1, srs=r^{-1} \rangle$  has an epimorphism onto 
$\langle\left\{e_{2}e_{4}\rho_{12}, (x_{3},x_{4})e_{3})\right\}\rangle$ and $\langle\left\{e_{2}e_{4}\rho_{12}, (x_{4},x_{5})e_{4})\right\}\rangle$
and hence these groups are  finite.

\end{pf}

\section{Some facts about \CAT\ spaces}

In this section we briefly present the main definitions and properties concerning \CAT\ metric spaces.
A detailed description of these spaces and their geometry can be found in  \cite{Haefliger}.

We start by reviewing the concept of geodesic spaces. Let $(X,d)$ be a metric space and $x, y\in X$. A {\bf geodesic} joining $x$ and $y$ is a map $c_{xy}:[0,l]\rightarrow X$, such that $c_{xy}(0)=x$, $c_{xy}(l)=y$ and $d(c_{xy}(t),c_{xy}(t'))=|t-t'|$ for all $t,t'\in[0,l]$. The image of~$c_{xy}$, denoted by $[x,y]$, is called a {\bf geodesic segment}. 
A metric space $(X,d)$ is said to be a {\bf geodesic space} if every two points in $X$ are joined by a geodesic. 
We say that $X$ is {\bf uniquely geodesic} if there is exactly one geodesic joining $x$ and $y$ for all $x,y\in X$. 

A {\bf geodesic ray} in $X$ is a map $c:[0, \infty)\rightarrow X$, such that $d(c(t),c(t'))=|t-t'|$ for all $t,t'\geq0$. A {\bf geodesic line} in $X$ is a map $c:\mathbb{R}\rightarrow X$ such that $d(c(t), c(t'))=|t-t'|$ for all $t,t'\in\mathbb{R}$.

A {\bf geodesic triangle} in $X$ consists of three points $p_{1},p_{2},p_{3}$ in $X$ and a choice of three geodesic segments $[p_{1}, p_{2}]$, $[p_{2}, p_{3}]$, $[p_{3}, p_{1}]$. Such a geodesic triangle will be denoted by $\Delta(p_{1}, p_{2}, p_{3})$. A triangle $\overline{\Delta}(\overline{p_{1}}, \overline{p_{2}}, \overline{p_{3}})$ in Euclidian space $\mathbb{R}^{2}$ is called a {\bf comparison triangle} for $\Delta(p_{1}, p_{2},$ $p_{3})$ if it is a geodesic triangle in $\mathbb{R}^{2}$ and if 
$d(p_{i}, p_{j})=d(\overline{p_{i}}, \overline{p_{j}})$ for $i, j=1,2,3$. A point $\overline{x}$ in $[\overline{p_{i}}, \overline{p_{j}}]$  is called a comparison point for $x\in[p_{i}, p_{j}]$ if  
$d(x, p_{i})=d(\overline{x}, \overline{p}_{i})$ and $d(x, p_{j})=d(\overline{x}, \overline{p}_{j})$. A geodesic triangle in $X$ is said to satisfy the $\mathbf{CAT(0)}$ {\bf inequality} if for all $x$ and $y$ in the geodesic triangle and all comparison points $\overline{x}$ and $\overline{y}$, the inequality $d(x,y)\leq d(\overline{x},\overline{y})$ holds. 

\begin{definition}
A metric space $X$ is called a $\mathbf{CAT(0)}$ \textbf{space} if $X$ is a geodesic space and all of its geodesic triangles satisfy the \CAT\ inequality. 
\end{definition}

One can easily verify from the definition of a \CAT\ space that these spaces are uniquely geodesic, therefore we use the notation $[x,y]$ for the geodesic segment between $x$ and $y$ in the \CAT\ space $X$. A subset $Y$ of a \CAT\ space $X$ is called {\bf convex} if for all $x$ and $y$ in $Y$ the geodesic segment $[x,y]$ is contained in $Y$. Indeed, convex subspaces of a \CAT\ space are obviously \CAT\ spaces. The {\bf diameter} of $Y$ is defined as ${\rm diam}\text{$(Y)=\sup\left\{d(x,y)\mid x,y\in Y\right\}$}$. The subset $Y$ is called {\bf bounded} if ${\rm diam}(Y)$ is finite. We also note, that the metric on a \CAT\ metric space is {\bf convex}, meaning that
for each pair of geodesics $c_1:[0, a_1]\rightarrow X$ and $c_2:[0, a_2]\rightarrow X$ with $c_1(0)=c_2(0)$ the inequality $d(c_1(ta_1), c_2(ta_2))\leq td(c_1(a_1), c_2(a_2))$ holds for all $t\in [0,1]$. 

The class of \CAT\ spaces is large. Perhaps the easiest examples of \CAT\ spaces besides $d$-dimensional Euclidean spaces $\mathbb{R}^d$ are metric trees, in particular simplicial trees by metrizing each edge of a simplicial tree as interval with length one. 

Let us mentioned an important property of \CAT\ spaces which will be needed later.
\begin{proposition}(\cite[\Romannum{2} 1.4]{Haefliger})
\label{Acyclic}
Any \CAT\ metric space is contractible, in particular all of its higher singular homology groups are trivial.
\end{proposition}

Now that we have introduced a class of spaces, we need, as in other mathematical theories, structure preserving maps. For a metric space $(X,d)$ an isometry $f: X\rightarrow X$ is a bijection such that $d(f(x), f(y))=d(x,y)$ for all $x$ and $y$ in $X$. The group of all isometries of $X$ will be denoted by \Isom. 
One easily checks that the fixed point set of an isometry of a \CAT\ space is closed and convex (or empty).

The following version of the Bruhat-Tits Fixed Point Theorem \cite[\Romannum{2} 2.8]{Haefliger} is crucial for our arguments. 
\begin{proposition}
\label{boundedOrbit}
 Let $G$ be a group acting on a complete \CAT\ space $X$ by isometries. Then the following conditions are equivalent: 
\begin{enumerate}
 \item[$(i)$] The group $G$ has a global fixed point.
 \item[$(ii)$] Each orbit of $G$ is bounded.
 \item[$(iii)$] The group $G$ has a bounded orbit.
\end{enumerate}
If the group $G$ satisfies one of the conditions above, then $G$ is called \textbf{ bounded on} $\mathbf{X}$.\\
\end{proposition}
The implications $i)\Rightarrow ii)$ and $ii)\Rightarrow iii)$ are trivial, and $iii)\Rightarrow i)$ is proven in \cite[\Romannum{2} 2.8]{Haefliger}.
$ $\\

The following corollary  is standard consequence of Proposition \ref{boundedOrbit}.

\begin{corollary}
\label{comm}
Let $G_{1}$, $G_{2}$ be groups, $X$ a complete \CAT\ space and 
\begin{align*}
\phi_{1}&:G_{1}\rightarrow\isom,\\
\phi_{2}&:G_{2}\rightarrow\isom 
\end{align*}
 homomorphisms. If $G_{1}$ and $G_{2}$ are bounded on $X$ and $\phi_{1}(g_{1})\circ\phi_{2}(g_{2})=\phi_{2}(g_{2})\circ\phi_{1}(g_{1})$ for all $g_{1}$ in $G_{1}$ and $g_{2}$ in $G_{2}$, then the map
\begin{align*}
 \phi_{1}\times\phi_{2}:G_{1}\times G_{2}&\rightarrow\isom\\
                           (g_{1}, g_{2})&\mapsto\phi(g_{1})\circ\phi(g_{2})
\end{align*}
is a homomorphism and $G_{1}\times G_{2}$ is bounded on $X$.
\end{corollary}

\section{Helly's Theorem for complete \CAT\ spaces and homological properties of nerves}

The purpose of this section is to verify one important result of convexity theory, Helly's Theorem, for the class of finite dimensional~\CAT\ spaces. 

\begin{HellyRTheorem}(\cite{Helly})
\label{ClasHelly}
Let  $\mathcal{S}$ be a finite family of non-empty closed convex subspaces of $\R^{d}$. If the intersection of each $(d+1)$-elements of $\mathcal{S}$ is non-empty, then $\bigcap\mathcal{S}$ is non-empty.
\end{HellyRTheorem}

There exist numerous different versions of this theorem for \CAT\ spaces in the literature, e.g. for finite families of non-empty convex open resp. closed subspaces, see \cite[3.2]{MappingClassFA}, 
\cite[2]{Debrunner}, \cite[3.2]{Farb} and \cite[5.3]{Kleiner}. 

We will study the proof by \textsc{Debrunner} for this result \cite{Debrunner}. \textsc{Debrunner} formulated and proved Helly's Theorem for a family of convex open subspaces of~$\R^{d}$. As we will see in this section, the same line of arguments as in \cite{Debrunner} also works for a family of convex open subspaces of a 
$d$-dimensional \CAT\ space. \textsc{Farb} observed in~\cite{Farb} that Helly's Theorem for  open convex subspaces of a $d$-dimensional \CAT\ space implies the version  for closed convex subspaces. Here we include a complete proof for this version.

We require the following definition. For a topological space $X$ we consider the reduced singular homology groups $\widetilde{\H}_{q}(X)$ for~$q\in\mathbb{Z}$. 

\begin{definition}
A topological space $X$ is said to be \textbf{acyclic} if~$\widetilde{\H}_{q}(X)=0$ for all~$q\in\mathbb{Z}$.
\end{definition}

In particular, if $X$ is acyclic then $X$ is non-empty and connected. For example, every contractible space is acyclic. Hence \CAT\ spaces are acyclic, see Proposition~\ref{Acyclic}.
$ $\\

The proof by \textsc{Debrunner} is based on the following proposition.

\begin{proposition}(\cite[Lemma $A_m$]{Debrunner})
\label{intersection}
 Let $X$ be a topological space and $\mathcal{S}$,~$|\mathcal{S}|\geq 2$ a finite family of open non-empty subspaces such that $\bigcap\mathcal{T}$ is acyclic for all $\mathcal{T}\subset\mathcal{S}$ with $|\mathcal{T}|=1,\ldots,|\mathcal{S}|-1$.
\begin{enumerate}
 \item[$(i)$] If $\bigcap\mathcal{S}$ is empty, then $\widetilde{\H}_{|\mathcal{S}|-2}(\bigcup\mathcal{S})\neq 0$.
 \item[$(ii)$] If $\bigcap\mathcal{S}$ is non-empty, then 
$\widetilde{\H}_{*}(\bigcup\mathcal{S})\cong\widetilde{\H}_{*-|\mathcal{S}|+1}(\bigcap\mathcal{S})$. In particular,
$\bigcup\mathcal{S}$ is acyclic if and only if $\bigcap\mathcal{S}$ is acyclic. 
\end{enumerate}
\end{proposition}
\begin{pf}
We prove both statements by induction on $m:=|\mathcal{S}|$.
Let $\mathcal{S}$ be equal to~$\left\{X_{1},X_{2}\right\}$. 
If~$\bigcap\mathcal{S}$ is empty, then $\bigcup\mathcal{S}$ is not connected and we have $\widetilde{\H}_{0}(\bigcup\mathcal{S})\neq 0$.
If~$\bigcap\mathcal{S}$ is non-empty, then we consider the reduced Mayer-Vietoris sequence for the pair $(X_{1}, X_{2})$.
\[
\ldots\rightarrow\widetilde{\H}_{q}(X_{1})\oplus\widetilde{\H}_{q}(X_{2})\rightarrow \widetilde{\H}_{q}(X_{1}\cup X_{2})\rightarrow\widetilde{\H}_{q-1}(X_{1}\cap X_{2})\rightarrow\widetilde{\H}_{q-1}(X_{1})\oplus\widetilde{\H}_{q-1}(X_{2})\rightarrow\ldots
\]
We know that $X_{1}$ and $X_{2}$ are acyclic, so we obtain
\[
 \ldots\rightarrow 0\rightarrow\widetilde{\H}_{q}(X_{1}\cup X_{2})\rightarrow\widetilde{\H}_{q-1}X_{1}\cap X_{2})\rightarrow 0\rightarrow\ldots 
\]
and therefore 
\[
 \widetilde{\H}_{q}(X_{1}\cup X_{2})\cong\widetilde{\H}_{q-1}(X_{1}\cap X_{2}).
\]
Now assume that $m>2$. Let $\mathcal{S}=\left\{X_{1},\ldots,X_{m}\right\}$ be a family of open subspaces such that the intersection of each $r$ members of this family is acyclic whenever~$r=1,\ldots,m-1$. We define
$U_{1}:=X_{1}\cup\ldots\cup X_{m-1}$, $U_{2}:=X_{m}$ and consider the reduced Mayer-Vietoris sequence for the pair $(U_{1}, U_{2})$
\[
\ldots\rightarrow\widetilde{\H}_{q}(U_{1})\oplus\widetilde{\H}_{q}(U_{2})\rightarrow 
\widetilde{\H}_{q}(U_{1}\cup U_{2})\rightarrow\widetilde{\H}_{q-1}(U_{1}\cap U_{2})\rightarrow\widetilde{\H}_{q-1}(U_{1})\oplus\widetilde{\H}_{q-1}(U_{2})\rightarrow\ldots
\]
The subspace $U_{2}$ is acyclic by assumption, and  $U_{1}$ is acyclic by part $(ii)$ of the induction hypothesis.
We have
\[
\ldots\rightarrow 0\rightarrow 
\widetilde{\H}_{q}(U_{1}\cup U_{2})\rightarrow\widetilde{\H}_{q-1}(U_{1}\cap U_{2})\rightarrow 0\rightarrow\ldots
\]
and therefore
\[
\widetilde{\H}_{q}(U_{1}\cup U_{2})\cong\widetilde{\H}_{q-1}(U_{1}\cap U_{2}).
\]
Now we define $\mathcal{S}'=\left\{X_{1}\cap X_{m},X_{2}\cap X_{m},\ldots, X_{m-1}\cap X_{m}\right\}$. This is a finite family of open subspaces such that the intersection of each $r$ members of this family is acyclic whenever $r=1,\ldots,m-2$. 

If $\bigcap\mathcal{S}=\bigcap\mathcal{S}'$ is empty, then we have 
\begin{eqnarray*}
\widetilde{\H}_{m-2}(\bigcup\mathcal{S})&=&\widetilde{\H}_{m-2}(U_1\cup U_2)\\
								 &\cong&\widetilde{\H}_{m-3}(U_{1}\cap U_{2})\\
							&\overset{U_1\cap U_2=\bigcup\mathcal{S}' }{=}&\widetilde{\H}_{m-3}(\bigcup\mathcal{S}')\\
						&\overset{\text{Ind. hyp. } (i)}{\neq} &0.
\end{eqnarray*}

If $\bigcap\mathcal{S}=\bigcap\mathcal{S}'$ is non-empty, then 
\begin{eqnarray*}
 \widetilde{\H}_{q}(\bigcup\mathcal{S})&=&\widetilde{\H}_q(U_1\cup U_2)\\
							&\cong&\widetilde{\H}_{q-1}(U_1\cap U_2)\\
							&\overset{U_1\cap U_2=\bigcup\mathcal{S}' }{=}&\widetilde{\H}_{q-1}(\bigcup\mathcal{S}')\\
							&\overset{\text{Ind. hyp. }(ii)}\cong&\widetilde{\H}_{q-1-(m-1)+1}(\bigcap\mathcal{S}')\\
							&\overset{\bigcap\mathcal{S}'=\bigcap\mathcal{S}}{=}&\widetilde{\H}_{q-m+1}(\bigcap\mathcal{S}).\end{eqnarray*}\end{pf}
Proposition \ref{intersection} gives the following topological version of Helly's Theorem.
\begin{theorem}(compare \cite[Thm. 2]{Debrunner})
\label{HomHelly}
Let $X$ be a topological space, $d$ a natural number and  $\mathcal{S}$ a finite family of open non-empty subspaces with the properties
\begin{enumerate}
 \item[$(i)$] $\widetilde{\H}_{q}(\bigcup\mathcal{T})=0$ for all $q\geq d$ and all $\mathcal{T}\subseteq\mathcal{S}$,
 \item[$(ii)$] $\bigcap\mathcal{T}$ is acyclic for $\mathcal{T}\subseteq\mathcal{S}$ with $|\mathcal{T}|=1,\ldots,d+1$, 
\end{enumerate}
then $\bigcap\mathcal{S}$ is acyclic.
\end{theorem}
\begin{pf}
Assume there exist families satisfying the hypotheses but not the conclusion. Let $\left\{X_{1},\ldots, X_{m}\right\}$ be such a family of minimal order. Using $(ii)$ we have $m\geq d+2$. This family satisfies hypothesis of Proposition \ref{intersection} by minimality of $m$.

If $X_{1}\cap\ldots\cap X_{m}$ is empty, then by Proposition \ref{intersection} $(i)$ we have $\widetilde{\H}_{m-2}(\bigcup\mathcal{S})\neq 0$. This contradicts $(i)$.

If $X_{1}\cap\ldots\cap X_{m}$ is non-empty, then there exists $q\geq 0$ with $\widetilde{\H}_{q}(X_{1}\cap\ldots\cap X_{m})\neq 0$. By Proposition \ref{intersection} $(ii)$ follows that 
$\widetilde{\H}_{q+m-1}(X_{1}\cup\ldots\cup X_{m})\neq 0$. We have $q+m-1\geq d$. This contradicts $(i)$.
\end{pf}

Using Theorem \ref{HomHelly} we can easily prove  the following version of Helly's Theorem for a family of convex open subspaces of a 
$d$-dimensional \CAT\ space. Recall, that here by dimension we mean the  covering dimension of a metric space. In contrast to that, the  {\bf compact dimension} of a space $X$ is defined as
\[
{\rm cdim}(X) = {\rm max}\left\{{\rm dim}(Y)\mid Y\subseteq X \text{ compact}\right\}. 
\]
We note that for a metric space $X$ we have clearly
\[
{\rm cdim}(X)\leq{\rm dim}(X),
\]
since ${\rm dim}(Y)\leq{\rm dim}(X)$ holds for all compact subsets. 

Before we turn to the proof we need the following result.
\begin{proposition}
\label{dimension}
Let $X$ be a $d$-dimensional \CAT\ space. Then the reduced singular homology groups $\widetilde{\rm H}_{q}(U)=0$ for all open subspaces $U\subseteq X$ and all $q\geq d$.
\end{proposition}
\begin{pf}
For $d=0$ there is nothing to prove. Assume that $d\geq 1$.
As shown by \textsc{Kleiner} in \cite[Thm. A]{Kleiner} one has for any \CAT\ space $X$
\[
{\rm cdim}(X)={\rm max}\left\{k\mid {\rm H}_k(U, V)\neq 0\text{ for some open pair } (U, V)\text{ in }X\right\}.
\]
We have $d={\rm dim}(X)\geq {\rm cdim}(X)$,
therefore we obtain $\widetilde{{\rm H}}_q(U,V)=0$ for all open pairs~$(U, V)$ in~$X$ and all~$q>d$.
Let $U\subseteq X$ be an open subspace. We consider the long exact sequence for the pair~$(X,U)$:
\[
\ldots\rightarrow \widetilde{{\rm H}}_{n+1}(X)\rightarrow \widetilde{{\rm H}}_{n+1}(X,U)\rightarrow\widetilde{{\rm H}}_{n}(U)\rightarrow\widetilde{{\rm H}}_{n}(X)\rightarrow\ldots
\]
The \CAT\ space $X$ is contractible, therefore $\widetilde{\H}_{*}(X)=0$ and we have
\[
\widetilde{{\rm H}}_{n+1}(X,U)\cong\widetilde{{\rm H}}_{n}(U).
\]
Using $\widetilde{{\rm H}}_q(X,U)=0$ for all $q>d$ we obtain $\widetilde{{\rm H}}_{q}(U)=0$ for all $q\geq d$.
\end{pf}

\begin{HellyOpenTheorem}
\label{HomHelly2}
Let~$X$ be a $d$-dimensional complete \CAT\ space and $\mathcal{S}$ a finite family of non-empty open convex subspaces. If the intersection of each $(d+1)$-elements of $\mathcal{S}$ is non-empty, then~$\bigcap\mathcal{S}$ is non-empty.
\end{HellyOpenTheorem}
\begin{proof}
The \CAT\ space $X$ has covering dimension $d$, therefore by Proposition~\ref{dimension},  we have $\widetilde{\H}_{q}(\bigcup\mathcal{T})=0$ for all $\mathcal{T}\subseteq\mathcal{S}$ and  $q\geq d$. Since the intersection of convex sets in $\mathcal{T}\subseteq\mathcal{S}$ with $|\mathcal{T}|=1,\ldots,d+1$ is non-empty and  convex,  $\bigcap\mathcal{T}$ is by Proposition \ref{Acyclic} contractible and hence acyclic. By Theorem \ref{HomHelly} it follows that $\bigcap\mathcal{S}$ is acyclic, in particular non-empty.
\end{proof}

For our application we require a variation of Helly's Theorem for closed convex subspaces of a $d$-dimensional \CAT\ space and we include a complete proof of this result here. Let us outline the structure of our proof: first we replace each of the closed convex subspaces by a bounded closed convex subspace. For this new family we then construct a swelling consisting of open convex bounded subspaces. Applying Helly's Theorem \ref{HomHelly2} for this constructed family we obtain a non-empty intersection of it, hence the intersection of a family we started with is also non-empty. 

We need the following definition.
\begin{definition}
Let $X$ be a topological space. A \textbf{swelling} of a family $(A_{i})_{i\in I}$ with $A_i\subseteq X$
is a family $(B_{i})_{i\in I}$ with $B_i\subseteq X$, such that $A_{i}\subseteq B_{i}$ for every $i\in I$ and for every finite subset $J\subseteq I$ we have
\[
 \bigcap\limits_{j\in J} A_{j}=\emptyset\text{ if and only if } \bigcap\limits_{j\in J} B_{j}=\emptyset.
\]
\end{definition}

Let us first recall an important property of a  \CAT\ space which we will need in the proof of the next proposition. By definition, a family of subsets $(A_i)_{i\in I}$ of a metric space is said to have the {\bf finite intersection property} if the intersection of each finite subfamily is non-empty. \textsc{Monod} proved in \cite[Thm. 14]{Monod} that a family consisting of bounded closed convex subsets of a complete \CAT\ space with the finite intersection property has a non-empty intersection.
\begin{proposition}
\label{distance}
 Let $X$ be a complete \CAT\ space and $A,B\subseteq X$ be non-empty closed convex subsets with $A\cap B=\emptyset$ and $A$ bounded, then 
\[
d(A,B):=\inf\left\{d(a,b)\mid a\in A, b\in B\right\}>0.
\]
\end{proposition}
\begin{pf}
 We assume that $d(A,B)=0$, then there exists a sequence $(a_{n})_{n\in\mathbb{N}}$ in~$A$ and a sequence $(b_{n})_{n\in\mathbb{N}}$ in $B$ such that $\lim\limits_{n} d(a_{n},b_{n})=0$. Let for all $n\in\mathbb{N}$, $A_{n}\subseteq A$ be the closed convex hull of the set $\left\{a_{k}\mid k\geq n\right\}$ and consider the family~$(A_{n})_{n\in\mathbb{N}}$. Since $A$ is bounded, this family consists of bounded closed convex subsets and has the finite intersection property, therefore the intersection of $\left\{A_{n}\mid n\in\mathbb{N}\right\}$ is non-empty. Using this fact we have $\bigcap\limits_{n\in\mathbb{N}}A_{n}\neq\emptyset$. Let $x$ be in $\bigcap\limits_{n\in\mathbb{N}}A_{n}\subseteq A$.
Further we define 
\[
B_{n}:=\left\{y\in X\mid d(y,B)\leq d(a_{n},b_{n})\right\}.
\]
The set $B_{n}$ is a closed convex set and $A_{n}\subseteq B_{n}$ for $n\in\mathbb{N}$. Therefore $x$ is in
$\bigcap\limits_{n\in\mathbb{N}}B_{n}$. Using $d(x,B)\leq d(a_{n},b_{n})$ for all $n\in\mathbb{N}$ we obtain $d(x,B)=0$. The set $B$ is closed, hence $x\in B$ and therefore $x\in A\cap B$. A contradiction.
\end{pf}

Using Proposition \ref{distance} we can now construct a swelling of a finite family of closed bounded convex subsets of a complete \CAT\ space which consists of open bounded convex subsets.

\begin{proposition}
\label{ClosedIntersection}
 Let $X$ be a complete \CAT\ space and $\mathcal{F}=\left\{F_{1},\ldots,F_{k}\right\}$ a finite family of non-empty closed bounded convex subsets. Then there exists a swelling  
$\mathcal{U}=\left\{U_{1},\ldots, U_{k}\right\}$ of $\mathcal{F}$ consisting of non-empty open bounded convex subsets.
\end{proposition}
\begin{pf}
We define 
\[
\mathcal{F}_{1}:=\left\{\bigcap\mathcal{G}\mid\mathcal{G}\subseteq\mathcal{F}, \bigcap\mathcal{G}\neq\emptyset, \bigcap\mathcal{G}\cap F_{1}=\emptyset\right\}.
\]
By Proposition \ref{distance}, $\min\left\{1,d(F_{1},S_{i})\mid S_{i}\in\mathcal{F}_{1}\right\}=\epsilon_{1}>0$. The family
\[
\mathcal{V}_{1}=\left\{V_{1},F_{2},\ldots,F_{k}\right\}
\]
where $V_{1}:=\left\{x\in X\mid d(x,F_{1})\leq \frac{\epsilon_{1}}{2}\right\}$ is a swelling of $\mathcal{F}$ which consists of non-empty closed bounded convex subsets. More precisely, the subset $V_{1}$ is convex because the subset $F_{1}$ and the \CAT\ metric are convex.

Now we assume that for $i\in\left\{1,\ldots,j\right\}$ the family $\mathcal{V}_{j}:=\left\{V_{1},\ldots,V_{j},F_{j+1},\ldots, F_{k}\right\}$ is defined and is a swelling of $\mathcal{F}$ which consists of non-empty closed bounded subsets.
We define 
\[
\mathcal{F}_{j+1}:=\left\{\bigcap\mathcal{G}\mid\mathcal{G}\subseteq\mathcal{V}_{j}, \bigcap\mathcal{G}\neq\emptyset, \bigcap\mathcal{G}\cap F_{j+1}=\emptyset\right\}.
\]
By  Proposition \ref{distance}, $\min\left\{1,d(F_{j+1},S_{i})\mid S_{i}\in\mathcal{F}_{j+1}\right\}=\epsilon_{j+1}>0$. The family
\[
\mathcal{V}_{j+1}=\left\{V_{1},\ldots,V_{j+1},F_{j+2},\ldots,F_{k}\right\}
\]
where 
$V_{j+1}:=\left\{x\in X\mid d(x,F_{j+1})\leq \frac{\epsilon_{j+1}}{2}\right\}$ is a swelling of $\mathcal{F}$ which consists of non-empty closed bounded convex subsets.
Thus, we can assume that $\mathcal{V}_{k}=\left\{V_{1},\ldots,V_{k}\right\}$ with $V_{i}=\left\{x\in X\mid d(x,F_{i})\leq\epsilon_{i}\right\}$ is defined and is a swelling of $\mathcal{F}$. The family 
$\mathcal{U}:=\left\{U_{1},\ldots, U_{k}\right\}$ with $U_{i}=\left\{x\in X\mid d(x,F_{i})<\frac{\epsilon_{i}}{2}\right\}$ is a swelling of $\mathcal{F}$ which consists of non-empty bounded open convex subsets.
\end{pf}

We are now ready to prove Helly's Theorem for a finite family of closed convex subspaces of a \CAT\ space.

\begin{HellyClosedTheorem}
\label{Helly}
 Let $X$ be a $d$-dimensional complete \CAT\ space and $\mathcal{S}$ a finite family of non-empty closed convex subspaces. If the intersection of each $(d+1)$-elements of $\mathcal{S}$ is non-empty, then $\bigcap\mathcal{S}$ is non-empty.
\end{HellyClosedTheorem}
\begin{pf}
For each subset $\mathcal{T}$ of $\mathcal{S}$ of order equal to $d+1$ we choose an element $p$ in $\bigcap\mathcal{T}$. Let the union of these elements be the set $\mathcal{P}$. This set is finite and we define 
\[
\mathcal{S}'=\left\{\overline{{\rm conv}}\left\{\mathcal{P}\cap S\right\}\mid S\in\mathcal{S}\right\},
\]
 where $\overline{{\rm conv}}\left\{\mathcal{P}\cap S\right\}$ is the closure of the convex hull of $\left\{\mathcal{P}\cap S\right\}$. The set $\mathcal{S}'$ consists of non-empty, closed bounded convex subspaces and the intersection of each $(d+1)$-elements of $\mathcal{S}'$ is non-empty. By Proposition \ref{ClosedIntersection} there exists a swelling $\mathcal{U}$ of $\mathcal{S}'$ which consists of non-empty open convex subspaces. By Helly's Theorem \ref{HomHelly2} it follows that $\bigcap \mathcal{U}$ is 
non-empty and  therefore $\bigcap\mathcal{S}'$ is non-empty. We have $\emptyset\neq\bigcap \mathcal{S}'\subseteq\bigcap \mathcal{S}$. This completes the proof.
\end{pf}
Our main technique in the proofs of Theorems A and B is based on the following crucial corollary. Indeed, it was \textsc{Farb} who discovered the connection between Helly's Theorem and the combinatorics of generating sets for a large class of groups.
\begin{Farb}
\label{HellyGroup}
Let $G$ be a group, $Y$ a finite generating set of $G$ and $X$ a complete $d$-dimensional \CAT\ space. Let
$\Phi:G\rightarrow\isom$ be a homomorphism. If each $(d+1)$-element subset of $Y$ has a fixed point in $X$, then $G$ has a fixed point in $X$.
\end{Farb}
\begin{pf}
Recall that $\Fix(y)\subseteq X$, the fixed point set of $y\in Y$, is a closed convex subset. Let $y_{1}$ and  $y_{2}$ be in $Y$, then $\Fix(y_{1})\cap\Fix({y_{2}})$ is equal to
$\Fix(\langle y_{1}, y_{2}\rangle)$ and therefore the statement immediately follows from Helly's Theorem for closed convex subspaces of a \CAT\ space \ref{Helly}.
\end{pf}

\section{Some facts about simplicial complexes and nerves}

In the previous section  we presented an important tool concerning global fixed point properties for isometric actions of groups, namely Farb's Fixed Point Criterion \ref{HellyGroup}. Using this criterion it remains to find a 'nice' generating set~$Y$ of~\Aut\ such that each of its $(d+1)$-element subsets has a fixed point. The purpose of this section is to present techniques to show that an infinite subgroup which is generated by $(d+1)$-elements of $Y$ has a fixed point. The methods presented here are based on certain simplicial complexes. 

Let us recall some basic properties of abstract simplicial complexes, see \cite{Haefliger} for details. A {\bf simplicial complex} $\Delta$ with a non-empty vertex set $\mathcal{V}$ is a collection of finite subsets of $\mathcal{V}$, called simplices, such that every one element subset of $\mathcal{V}$ is a simplex and $\Delta$ is closed under taking subsets. Let $A\in\Delta$ be a simplex. The cardinality $r$ of $A$ is called the {\bf rank} of $A$ and $r-1$ is called the {\bf dimension} of $A$. The {\bf dimension} of $\Delta$ is defined as: 
$\dim(\Delta):=\sup\left\{\dim(A)\mid A\in\Delta\right\}$. As usual, we denote by $|\Delta|$ the geometric realization of a simplicial complex $\Delta$.

In the following we need one basic construction that allow us to produce new simplicial complexes from old ones.
\begin{definition}
Let $K_{1}$, $K_{2}$ be simplicial complexes with vertex sets $\mathcal{V}_{1}$, $\mathcal{V}_{2}$. The \textbf{join} $K_{1}*K_{2}$ of $K_{1}$ and $K_{2}$ is a simplicial complex with vertex set equal to the union $\mathcal{V}_{1}\cup\mathcal{V}_{2}$ and $A\subseteq\mathcal{V}_{1}\cup\mathcal{V}_{2}$ is a simplex in $K_{1}*K_{2}$ if and only if $A=A_{1}\cup A_{2}$ where $A_{1}$ is a simplex in $K_{1}$ and $A_{2}$ is a simplex in $K_{2}$. 
\end{definition}

For example, the join of a standard $n$-simplex with a vertex set $\left\{0, 1, \ldots, n\right\}$, denoted by $\Delta_n$, and a standard  $m$-simplex with a vertex set 
$\left\{n+1, \ldots, n+m+1\right\}$ is a $(n+m+1)$-simplex with a vertex set $\left\{0, 1, \ldots, n+m+1\right\}$, see for example Figure 1.

\begin{figure}[h]
\begin{center}
\begin{tikzpicture}
\draw[line width=1.5pt]  (0,0)--(1,0);

\node at (0.5, -0.5) {\small{$|\Delta_1|$}};
\node at (0.5, 1.25) {\small{$|\Delta_0|$}};
\node at (2.5, -0.5) {\small{$|\Delta_1*\Delta_0|$}};

\node at (0,-0.3) {\tiny{$\left\{0\right\}$}};
\node at (1,-0.3) {\tiny{$\left\{1\right\}$}};
\node at (0.8,0.75) {\tiny{$\left\{2\right\}$}};
\filldraw[green] (2,0)--(3,0)--(2.5, 0.75);
\node at (0.5,0.2) {\tiny{$\left\{0,1\right\}$}};
\node at (3.2,0.4) {\tiny{$\left\{1,2\right\}$}};

\draw[line width=1.5pt]  (2,0)--(3,0);
\draw[line width=1.5pt]  (2,0)--(2.5,0.75);
\draw[line width=1.5pt]  (3,0)--(2.5,0.75);

\draw[fill=black]  (0, 0) circle (2pt);
\draw[fill=black]  (1, 0) circle (2pt);
\draw[fill=black]  (0.5, 0.75) circle (2pt);
\draw[fill=black]  (2, 0) circle (2pt);
\draw[fill=black]  (3, 0) circle (2pt);
\draw[fill=black]  (2.5, 0.75) circle (2pt);
\end{tikzpicture}
\end{center}
\caption{}
\end{figure}

If the geometric realisation of a simplicial complex $K_i$ is homeomorphic to a  sphere of dimension $d_i$ for $i=1,2$, then the geometric realisation of $K_1*K_2$ is homeomorphic to a sphere of dimension $d_1+d_2-1$. In particular, we have:
\[
|\partial\Delta_n*\partial\Delta_m|\cong\mathbb{S}^{n+m-1}
\]
for $n,m>0$, where  $\partial\Delta_n$ resp. $\partial\Delta_m$ is a boundary of $\Delta_n$ resp. $\Delta_m$, see for example Figure 2.

\begin{figure}[h]
\begin{center}
\begin{tikzpicture}

\node at (0, -0.5) {\small{$|\Delta_1|$}};
\node at (1, -0.5) {\small{$|\Delta_1|$}};

\node at (3.3, -0.5) {\small{$|\partial\Delta_1*\partial\Delta_1|\cong\mathbb{S}^1$}};

\draw[line width=1.5pt]  (0,0)--(1,1);

\draw[line width=1.5pt]  (1,0)--(0,1);

\draw[line width=1.5pt]  (2,0)--(3,0);
\draw[line width=1.5pt]  (3,0)--(3,1);
\draw[line width=1.5pt]  (3,1)--(2,1);
\draw[line width=1.5pt]  (2,1)--(2,0);

\draw[fill=green]  (0, 0) circle (2pt);
\draw[fill=red]  (1, 0) circle (2pt);
\draw[fill=red]  (0, 1) circle (2pt);
\draw[fill=green]  (1, 1) circle (2pt);
\draw[fill=green]  (2, 0) circle (2pt);
\draw[fill=red]  (3, 0) circle (2pt);
\draw[fill=red]  (2, 1) circle (2pt);
\draw[fill=green]  (3, 1) circle (2pt);
\end{tikzpicture}
\end{center}
\caption{}
\end{figure}

In the following we want to represent a family of subspaces of a topological space by a combinarial structure. For this reason we need the following definition.
\begin{definition}
 Let $X$ be a set and $\mathcal{S}$ a family of subsets of $X$. The \textbf{nerve} $\mathcal{N}(\mathcal{S})$ is the simplicial complex whose vertex set is $\mathcal{S}$ and whose non-empty simplices are all finite subsets $\left\{S_{1},\ldots,S_{k}\right\}\subseteq\mathcal{S}$ with $S_{1}\cap\ldots\cap S_{k}\neq\emptyset$. 
\end{definition}

For example, let $\mathcal{S}=\left\{S_0, \ldots, S_k\right\}$ be a set of non-empty closed convex subspaces of a $d$-dimensional complete \CAT\ space. We consider  $\mathcal{N}(\mathcal{S})$  as a subcomplex of the standard $k$-simplex $\Delta_k$. If the nerve $\mathcal{N}(\mathcal{S})$ contains the full $d$-skeleton of $\Delta_k$, then by Helly's Theorem \ref{Helly} it follows that 
$|\mathcal{N}(\mathcal{S})|\cong|\Delta_k|$. 

The main use of the nerve is the following proposition, due to \textsc{McCord}. Recall, by definition, a cover of a topological space is said to be {\bf point-finite} if every point of this space is contained in only finitely many sets of this cover. 
\begin{proposition}(\cite[2]{McCord})
\label{Nerve}
Let $Y$ be a topological space and \ $\mathcal{U}$ be a point-finite open
cover of $Y$ such that the intersection of any finite subcollection of  $\mathcal{U}$ is either empty
or contractible. Then $\H_{*}(\mathcal{N}(\mathcal{U}))\cong\H_{*}(Y)$.
\end{proposition}
The next important result is Theorem \ref{fix}, the proof relies on a
Proposition \ref{Nerve} and the following result.
\begin{proposition}(\cite[3.3]{MappingClassFA})
\label{JoinNerve}
 Let $X$ be a complete \CAT\ space and let $S_{1},\ldots, S_{l}$ be subsets of \Isom\ such that $[S_{i}, S_{j}]=1$ holds for all $1\leq i<j\leq l$.
Let $\mathcal{F}_{i}=\left\{\Fix(s)\mid s\in S_{i}\right\}$ and $\mathcal{N}_{i}=\mathcal{N}(\mathcal{F}_{i})$. Put 
$\mathcal{N}=\mathcal{N}(\mathcal{F}_{1}\cup\ldots\cup\mathcal{F}_{l})$. Then we have
\[
\mathcal{N}=\mathcal{N}_{1}*\ldots*\mathcal{N}_{l}.
\]
\end{proposition}
\begin{proof}
We first note that the vertex sets of $\mathcal{N}$ and $\mathcal{N}_{1}*\ldots*\mathcal{N}_{l}$ are equal.

Let $A=\left\{\Fix(s_{1}),\ldots,\Fix(s_{k})\right\}$ be a simplex in $\mathcal{N}$. We write the set $A$ as $A=A_{1}\cup\ldots\cup A_{l}$ with $A_{i}\subseteq\mathcal{F}_{i}$ for $i$ in $\left\{1,\ldots,l\right\}$. The intersection $\bigcap A_{i}$ is non-empty for all $i$ in $\left\{1,\ldots,l\right\}$ and therefore the set $A_{i}$ is a simplex in $\mathcal{N}_{i}$ for every $i$ in $\left\{1,\ldots,l\right\}$. It follows that the subset $A$ is a union of simplices in $\mathcal{N}_{i}$ and therefore a simplex in $\mathcal{N}_{1}*\ldots*\mathcal{N}_{l}$.
We have shown that $\mathcal{N}\subseteq\mathcal{N}_{1}*\ldots*\mathcal{N}_{l}$.

Now we prove the other inclusion. Let $B$ be a simplex in $\mathcal{N}_{1}*\ldots*\mathcal{N}_{l}$. We know that $B=B_{1}\cup\ldots\cup B_{l}$ where $B_{i}$ is a simplex in $\mathcal{N}_{i}$ for $i$ in $\left\{1,\ldots,l\right\}$. Now we have to show that $\bigcap B$ is non-empty. 
We consider the set $S'_{i}:=\left\{s\in S_{i}\mid\Fix(s)\in B_{i}\right\}$ and the group which is generated by $S'_{i}$ for $i$ in $\left\{1,\ldots,l\right\}$. The subset $B_{i}$ is a simplex in $\mathcal{N}_{i}$ and therefore the fixed point set $\Fix(\langle S'_{i}\rangle)$ is non-empty for all $i$ in $\left\{1,\ldots,l\right\}$. Next we note that $[\langle S'_{i}\rangle,\langle S'_{j}\rangle]=1$ for $i\neq j$. It follows from Corollary \ref{comm} that
\[
\bigcap^{l}\limits_{i=1}\Fix(\langle S'_{i}\rangle)=\Fix(\bigcup^{l}\limits_{i=1}(\langle S'_{i}\rangle))\neq\emptyset.
\]
In particular the set $\bigcap B$ is non-empty and therefore $B$ is a simplex in $\mathcal{N}$. This completes the proof.
\end{proof}

\begin{theorem}(\cite[3.4]{MappingClassFA})
\label{fix}
 Let $k_{1},\ldots, k_{l}$ be in $\mathbb{N}_{>0}$ and let $X$ be a $d$-dimensional complete \CAT\ space with $0<d<k_{1}+\ldots +k_{l}$. Let $S_{1},\ldots,S_{l}$  be subsets of \Isom\ such that $[S_{i}, S_{j}]=1$ for all $1\leq i<j\leq l$. If each $k_{i}$-element subset of $S_{i}$ has a fixed point in $X$ for all $i\in \left\{1,\ldots, l\right\}$, then for some $j\in\left\{1,\ldots, l\right\}$ every finite subset of $S_{j}$ has a fixed point.
\end{theorem}
\begin{proof}
Assume this is false, i.e. for each $i\in\left\{1,\ldots, l\right\}$ let $k^{'}_{i}\geq k_{i}$ be minimal such that there exists a $(k^{'}_{i}+1)$-element subset 
\[
T_{i}=\left\{s_{i,1},\ldots, s_{i,k^{'}_{i}+1}\right\}\subseteq S_{i}
\]
with empty fixed point set. By minimality of $k^{'}_i$, we know that each $k^{'}_{i}$-element subset of $T_{i}$ has a fixed point. Therefore the nerve of 
\[
\mathcal{F}_{i}=\left\{\Fix(s_{i,1}),\ldots,\Fix(s_{i,k^{'}_{i}+1})\right\}
\] 
is the boundary of a $k^{'}_{i}$-simplex. 
It follows from Proposition \ref{JoinNerve} that 
\[
\mathcal{N}(\mathcal{F}_{1}\cup\ldots\cup\mathcal{F}_{l})\cong\partial\Delta_{k'_{1}}*\ldots*\partial\Delta_{k'_{l}}.
\] 
The geometric realisation of this nerve is homeomorphic to a sphere, hence
\[
 |\mathcal{N}(\mathcal{F}_{1}\cup\ldots\cup\mathcal{F}_{l})|\cong\mathbb{S}^{k'_{1}+\ldots+k'_{l}-1}.
\]
Therefore the singular homology groups of the above spaces are isomorphic, i.e.
\[
 \H_{*}(|\mathcal{N}(\mathcal{F}_{1}\cup\ldots\cup\mathcal{F}_{l})|)\cong\H_{*}(\mathbb{S}^{k'_{1}+\ldots+k'_{l}-1}).
\]
Now we replace $\left\{\mathcal{F}_{1},\ldots,\mathcal{F}_{l}\right\}$ by a family consisting of a bounded convex closed subsets, as in the proof of Theorem \ref{Helly} and then by a swelling $\left\{\mathcal{F}'_{1},\ldots, \mathcal{F}'_{l}\right\}$ consisting of bounded convex open subsets, see Proposition \ref{ClosedIntersection}. We have
\[
\mathcal{N}(\mathcal{F}_{1}\cup\ldots\cup\mathcal{F}_{l})\cong\mathcal{N}(\mathcal{F}'_{1}\cup\ldots\cup\mathcal{F}'_{l}). 
\]
Using Proposition \ref{Nerve} we obtain  
\begin{align*}
\H_{k'_{1}+\ldots+k'_{l}-1}(\mathcal{F}'_{1}\cup\ldots\cup\mathcal{F}'_{l})&\cong\H_{k'_{1}+\ldots+k'_{l}-1}(|\mathcal{N}(\mathcal{F}'_{1}\cup\ldots\cup\mathcal{F}'_{l})|)\\
&\cong\H_{k'_{1}+\ldots+k'_{l}-1}(|\mathcal{N}(\mathcal{F}_{1}\cup\ldots\cup\mathcal{F}_{l})|)\\
&\cong\H_{k'_{1}+\ldots+k'_{l}-1}(\mathbb{S}^{k'_{1}+\ldots+k'_{l}-1}).
\end{align*}
Because the \CAT\ space $X$ is $d$-dimensional, we have by Proposition \ref{dimension} that the singular homology groups ${\rm H}_{q}(\mathcal{F}'_{1}\cup\ldots\cup\mathcal{F}'_{l})=0$ for  all $q\geq d$.
We have the inequality $k'_{1}+\ldots+k'_{l}-1\geq d$, in particular $\H_{k'_{1}+\ldots+k'_{l}-1}(\mathcal{F}'_{1}\cup\ldots\cup\mathcal{F}'_{l})\cong 0$.  This contradicts 
\[                                                                                                                                             
\H_{k'_{1}+\ldots+k'_{l}-1}(\mathbb{S}^{k'_{1}+\ldots+k'_{l}-1})\cong\Z.
\]
\end{proof}
The following consequence of Theorem \ref{fix}  is a crucial tool for proving global fixed point results for infinite subgroups.
\begin{corollary}(\cite[3.6]{MappingClassFA})
\label{conjugates}
 Let $k$ and $l$ be in $\mathbb{N}_{>0}$ and let $X$ be a complete $d$-dimensional \CAT\ space, with $d<k\cdot l$. Let $S$ be a subset of \Isom\ and let $S_{1},\ldots, S_{l}$ be conjugates of $S$ such that $[S_{i}, S_{j}]=1$ for $i\neq j$. If each $k$-element subset of $S$ has a fixed point in $X$, then each finite subset of $S$ has a fixed point in $X$.
\end{corollary}
\begin{proof}
This is clear from Theorem \ref{fix} since the fixed point sets of the sets $S_{i}$ are conjugate.
\end{proof}

\section{Proof of Theorem A}

Now we have all the ingredience to prove Theorem A. 
\begin{NewTheoremA}
\label{MainTheorem}
If $n\geq 4$ and $d<{\rm min}\left\{k\left\lfloor\frac{n}{k+2}\right\rfloor\mid k=2,\ldots, d+1\right\}$, then \Aut\ has property $\F\mathcal{A}_{d}$. In particular, if $n\geq 4$ and $d<2\left\lfloor\frac{n}{4}\right\rfloor-1$, then \Aut\ has property~$\F\mathcal{A}_{d}$. 
\end{NewTheoremA}
\begin{proof}
Let $X$ be a $d$-dimensional complete \CAT\ space and 
\[
 \Phi:\Autn\rightarrow\isom
\]
an action of \Aut\ on $X$. By Proposition \ref{GenAut} the group \Aut\ is generated by the set 
\[
Y_{2}:=\left\{(x_{1}, \ x_{2})e_{1}e_{2}, \ (x_{2},x_{3})e_{1},\ (x_{i}, x_{i+1}),\ e_{2}\rho_{12}, \ e_{n}\mid i=3,\ldots,n-1\right\}.
\]
Let us outline the structure of the proof: combining Bruhat-Tits Fixed Point Theo-  rem \ref{boundedOrbit} with Corollaries \ref{comm} and \ref{conjugates} we show the following: if $k\leq d+1$ and $d<{\rm min}\left\{k\left\lfloor\frac{n}{k+2}\right\rfloor\mid k=2,\ldots,d+1\right\}$, then each $k$-element subset of $Y_{2}$ has a fixed point. Then by Farb's Fixed point Criterion \ref{HellyGroup} the action $\Phi$ has a global fixed point.

As seen in Proposition \ref{GenAut} each element in $Y_{2}$ is an involution and the order of the product of two elements is finite. Let us consider the Coxeter group 
\[
W=\langle Y_{2}\mid (fg)^{\ord(fg)}=1, f,g\in Y_2\rangle
\]
whose Coxeter diagram looks as follows. 
\begin{figure}[h]
\begin{center}
\begin{tikzpicture}
 \draw (0,0.1)--(0.9,1);
 \draw (0,0.1)--(0,1.9);
 \draw (0,1.9)--(0.9,1);
 \draw (1.1,1)--(2.4,1);
 \draw (2.6,1)--(3.9,1);
 \draw (4.1,1)--(5.4,1);
 \draw (5.6,1)--(6.9,1);
 \draw (8.6,1)--(9.9,1);
 \draw (10.1,1)--(11.4,1);
 \draw (0,0) circle (3pt);
 \draw (0,-0.5) node {\footnotesize{$e_{2}\rho_{12}$}};
 \draw (0,2) circle (3pt);
 \draw (0,2.5) node {\footnotesize{$(x_{1},x_{2})e_{1}e_{2}$}};
 \draw (1,1) circle (3pt);
 \draw (1.25,0.5) node {\footnotesize{$(x_{2},x_{3})e_{1}$}};
 \draw (2.5,1) circle (3pt);
 \draw (2.5,1.5) node {\footnotesize{$(x_{3},x_{4})$}};
 \draw (4,1) circle (3pt);
 \draw (4,0.5) node {\footnotesize{$(x_{4},x_{5})$}};
 \draw (5.5,1) circle (3pt);
 \draw (5.5,1.5) node {\footnotesize{$(x_{5},x_{6})$}};
 \draw (7.80,1) node {$\ldots$};
 \draw (10,1) circle (3pt);
 \draw (10,0.5) node {\footnotesize{$(x_{n-1},x_{n})$}};
 \draw (11.5,1) circle (3pt);
 \draw (11.5,1.5) node {\footnotesize{$e_{n}$}}; 
 \draw (0.5,0.85) node {\footnotesize{$4$}};
 \draw (0.5,1.75) node {\footnotesize{$6$}};
 \draw (10.85,1.25) node {\footnotesize{$4$}};
 \draw (1.8,1.25) node {\footnotesize{$6$}};
\end{tikzpicture}
\caption{}
\end{center}
\end{figure}

In particular, we obtain an epimorphism
$\pi: W\rightarrow{\rm Aut}(F_n)$ and an action
\[
\Phi\circ\pi: W\rightarrow{\rm Aut}(F_n)\rightarrow{\rm Isom}(X).
\] 

It is obvious that if a subgroup of $W$ has a fixed point, then the image of this subgroup under $\pi$, a subgroup in \Aut, also has a fixed point.
For $k=2$ we know by Proposition \ref{GenAut} that each pair of the generating set $Y_2$ of \Aut\ generates a finite subgroup, therefore by the Bruhat-Tits Fixed Point Theorem~\ref{boundedOrbit} we obtain that each $2$-element subset of $Y_2$ has a fixed point. Now we assume that each $k$-element subset of $Y_2$ has a fixed point. Let $Y'$ be a $(k+1)$-element subset of $Y_{2}$.

If $e_{2}\rho_{12}$ is not in $Y'$, then it follows by Proposition \ref{GenAut} that $\langle Y'\rangle$ is a finite subgroup of \Aut\ and this subgroup has by Bruhat-Tits Fixed Point Theorem~\ref{boundedOrbit} a fixed point. 

If $e_{2}\rho_{12}$ is in $Y'$, we consider the Coxeter diagram of $\langle Y'\rangle\subseteq W$. If it is not connected, then it follows from hypothesis and from Corollary \ref{comm} that $\langle Y'\rangle $ has a fixed point. If the Coxeter diagram of $\langle Y'\rangle\subseteq W$ is connected, then we have the following cases:
\begin{enumerate}
 \item $Y'=\left\{e_{2}\rho_{12},(x_{1},x_{2})e_{1}e_{2},(x_{2}, x_{3})e_{1}, (x_{3},x_{4}),\ldots,(x_{k},x_{k+1})\right\}$,
 \item $Y'=\left\{e_{2}\rho_{12}, (x_{2},x_{3})e_{1}, (x_{3},x_{4}),\ldots,(x_{k+1},x_{k+2})\right\}$.
\end{enumerate}
The involution $e_{n}$ is not in $Y'$: assume that $e_{n}$ is contained in $Y'$, then $Y'$ consists of at least $n$ elements.  Therefore we must have $k+1\geq n$ which contradicts our assumption that
 $k+1\leq d+1<2\left\lfloor\frac{n}{4}\right\rfloor+1$. 

If $Y'$ is equal to $\left\{e_{2}\rho_{12},(x_{1},x_{2})e_{1}e_{2},(x_{2}, x_{3})e_{1},(x_{3},x_{4}),\ldots,(x_{k},x_{k+1})\right\}$, then  we define the permutations 
\[
 \tau_{i}:=(x_{1}, x_{(k+1)\cdot(i-1)+1})(x_{2}, x_{(k+1)\cdot(i-1)+2})\ldots(x_{k+1},x_{(k+1)\cdot(i-1)+k+1})
\]
and the sets 
\[
 S_{i}:=\tau_{i}Y'\tau^{-1}_{i}
\]
for $i\in\left\{1,\ldots,\left\lfloor\frac{n}{k+1}\right\rfloor\right\}$. The sets $S_{1},\ldots, S_{\left\lfloor\frac{n}{k+1}\right\rfloor}$ have the property that $[S_{i},S_{j}]=1$ for~$i\neq j$ as they act non-trivially only on disjoint subsets of $X$. By the assumption each $k$-element subset of $Y'$ has a fixed point and it follows from Corollary \ref{conjugates} that for $d<k\left\lfloor\frac{n}{k+1}\right\rfloor$ the set $Y'$ has a fixed point. 

If $Y'$ is equal to $\left\{e_{2}\rho_{12}, (x_{2},x_{3})e_{1}, (x_{3},x_{4}),\ldots,(x_{k+1},x_{k+2})\right\}$, then we define
the permutations 
\[
 \sigma_{i}:=(x_{1}, x_{(k+2)\cdot(i-1)+1})(x_{2}, x_{(k+2)\cdot(i-1)+2})\ldots(x_{k+2},\\x_{(k+2)\cdot(i-1)+k+2})
\]
and the sets 
\[
 T_{i}:=\sigma_{i}Y'\sigma^{-1}_{i}.
\]
for $i\in\left\{1,\ldots,\left\lfloor\frac{n}{k+2}\right\rfloor\right\}$. 
With similar arguments as above it follows that for $d<k\left\lfloor\frac{n}{k+2}\right\rfloor$ the set $Y'$ has a fixed point. 

So far we have shown that if $n\geq 4$ and $d< {\rm min}\left\{k\left\lfloor\frac{n}{k+2}\right\rfloor\mid k=2,\ldots,d+1\right\}$, then each $(d+1)$-element subset of $Y_2$ has a fixed point. By Farb's Fixed Point Criterion~\ref{HellyGroup} it follows that \Aut\ has a global fixed point.
An easy calculation shows: 
\[
2\left\lfloor\frac{n}{4}\right\rfloor-1\leq{\rm min}\left\{k\left\lfloor\frac{n}{k+2}\right\rfloor\mid k=2,\ldots,d+1\right\}.
\]
This completes the proof.
\end{proof}

Note that, as an immediate corollary of Theorem A, we obtain a similar result for~${\rm GL}_n(\mathbb{Z})$.
\begin{corollary}
If $n\geq 4$ and $d<{\rm min}\left\{k\left\lfloor\frac{n}{k+2}\right\rfloor\mid k=2,\ldots, d+1\right\}$, then  ${\rm GL}_n(\mathbb{Z})$ has property $\F\mathcal{A}_{d}$. In particular, if $n\geq 4$ and $d<2\left\lfloor\frac{n}{4}\right\rfloor-1$, then ${\rm GL}_n(\mathbb{Z})$ has property~$\F\mathcal{A}_{d}$. 
\end{corollary}

\section{Proof of Theorem B}
Using the result of Theorem A, we prove
\begin{NewTheoremB}
\label{FixedPointSAut}
If $n\geq 5$ and $d<{\rm min}\left\{k\left\lfloor\frac{n-1}{k+2}\right\rfloor\mid k=2,\ldots,d+1\right\}$, then \SAut\ has property $\F\mathcal{A}_{d}$.
In particular, if $n\geq 5$ and $d<2\left\lfloor\frac{n-1}{4}\right\rfloor-1$, then \SAut\ has pro- perty~$\F\mathcal{A}_{d}$.
\end{NewTheoremB}
\begin{pf}
Let $X$ be a $d$-dimensional complete \CAT\ space and 
\[
 \Phi:\SAutn\rightarrow\isom
\]
an action of \SAut\ on $X$. By Proposition \ref{finite2} the group \SAut\ is generated by the set
\[
Y_{4}:=\left\{(x_{1}, x_{2})e_{1}e_{2}e_{3},\  (x_{2},x_{3})e_{1},\ (x_{i}, x_{i+1})e_{i},\ e_{2}e_{4}\rho_{12}, \ e_{3}e_{4}\mid i=3,\ldots, n-1\right\}.
\]
If $n\leq 8$, then $d<2$ and the conclusion of Theorem B follows from Proposition \ref{finite2} and Farb's Fixed point Criterion \ref{HellyGroup}. We hence may assume that $n\geq 9$.
We show again the following: if $k\leq d+1$ and $d<{\rm min}\left\{k\left\lfloor\frac{n-1}{k+2}\right\rfloor\mid k=2,\ldots d+1\right\}$, then each $k$-element subset of $Y_{4}$ has a fixed point. 

Let us consider the Coxeter-like diagram for the set $Y_{4}$. We draw a graph with~$Y_{4}$ as vertex set, joining vertices $f$ and $g$ by an edge iff $[f,g]\neq 1$.
\begin{center}
\begin{figure}[h]
\begin{tikzpicture}
 \draw[fill=black] (3,-3) circle (3pt);
 \draw (3,-3.3) node {\footnotesize{$e_{3}e_{4}$}};
 \draw (1.6,0)--(3,-2.9);
 \draw (3,-2.9)--(4.4,0);
 \draw (3.1,0)--(4.4,0);
 \draw[fill=black] (1.5,0) circle (3pt);
 \draw (0.5,0) node {\footnotesize{$(x_{2},x_{3})e_{1}$}};
 \draw (1.6,0)--(2.9,0);
 \draw[fill=black] (3,0) circle (3pt);
 \draw (2.5,0.3) node {\footnotesize{$e_{2}e_{4}\rho_{12}$}};
 \draw[fill=black] (4.5,0) circle (3pt);
 \draw (4.5,-0.3) node {\footnotesize{$(x_{4},x_{5})e_{4}$}};
 \draw (4.6,0)--(5.9,0);
 \draw[fill=black] (6,0) circle (3pt);
 \draw (6,0.3) node {\footnotesize{$(x_{5},x_{6})e_{5}$}};
 \draw (6.1,0)--(7.4,0);
 \draw (8.3,0) node {$\ldots$};
 \draw (9.1,0)--(10.4,0);
 \draw[fill=black] (10.5,0) circle (3pt);
 \draw (10.5,-0.3) node {\footnotesize{$(x_{n-1},x_{n})e_{n-1}$}};
 \draw[fill=black] (3,1.5) circle (3pt);
 \draw (3,1.8) node {\footnotesize{$(x_{3},x_{4})e_{3}$}};
 \draw (1.6,0)--(2.9,1.5);
 \draw (3.1,1.5)--(4.4,0);
 \draw (3,0.1)--(3,1.4);
 \draw[fill=black] (3,-1.5) circle (3pt);
 \draw (3,-1.8) node {\footnotesize{$(x_{1},x_{2})e_{1}e_{2}e_{3}$}};
 \draw (1.6,0)--(2.9,-1.5);
 \draw (3,-0.1)--(3,-1.4);
\end{tikzpicture}
\caption{}
\end{figure}
\end{center}

If $k$ is equal to $2$, then we know by Proposition \ref{finite2} and Bruhat-Tits Fixed Point Theorem \ref{boundedOrbit} that each $2$-element subset of $Y_{4}$ has a fixed point.

Now we assume that each $k$-element subset of $Y_4$ has a fixed point. Let $Y'$ be a  $(k+1)$-element subset of $Y_{4}$. 
If $e_{2}e_{4}\rho_{12}$ is not in $Y'$, then it follows from Proposition~\ref{finite2} that $\langle Y'\rangle$ is a finite subgroup of \SAut\ and this subgroup has by Bruhat-Tits Fixed Point Theorem~\ref{boundedOrbit} a fixed point. If $e_{2}e_{4}\rho_{12}$ is in $Y'$ then we have the following cases:
\begin{enumerate}
\item If there exists a non-empty proper subset $Y''$ of $Y'$ with the property 
\[
[Y'',Y'-Y'']=1,
\]
then it follows by the assumption and by Corollary \ref{comm} that $Y'$ has a fixed point.
\item Otherwise we consider the determinant homomorphism
\[
 \det:\Autm\rightarrow\GL_{n-1}(\Z)\rightarrow\Z_{2}
\]
and we define
\[
 \Psi:\Autm\rightarrow\SAutn
\]
as follows
\[
 f\mapsto f'
\]
\[
\begin{matrix}
f'(x_{k}):=\begin{cases} f(x_{k}) & \mbox{if $k=1,\ldots,n-1$.}  \\ 
	          x^{\det(f)}_{k} & \mbox{if $k=n.$} 
\end{cases}
\end{matrix}
\]
The homomorphism $\Psi$ is injective and $Y'$ is contained in $\im(\Psi)$ because the element $(x_{n-1},x_{n})e_{n-1}$ is not contained in  $Y'$. More precisely, 
 assume that $(x_{n-1},x_{n})e_{n-1}\in Y'$, therefore $k+1\geq n-3$ which contradicts our assumption.

By Theorem A the group $\im(\Psi)\subseteq{\rm SAut}(F_n)$ has a global fixed point and therefore $Y'$ has a fixed point. 
\end{enumerate}
Again, by Farb's Fixed Point Criterion~\ref{HellyGroup} it follows that \SAut\ has a global fixed point. 
An easy calculation shows: 
\[
2\left\lfloor\frac{n-1}{4}\right\rfloor-1\leq {\rm min}\left\{k\left\lfloor\frac{n-1}{k+2}\right\rfloor\mid k=2,\ldots,d+1\right\}.
\]
This finishes the proof.
\end{pf}

Note that, as an immediate corollary of Theorem B we obtain a similar result  for~${\rm SL}_n(\mathbb{Z})$.
\begin{corollary}
If $n\geq 5$ and $d<{\rm min}\left\{k\left\lfloor\frac{n-1}{k+2}\right\rfloor\mid k=2,\ldots,d+1\right\}$, then ${\rm SL}_n(\mathbb{Z})$ has property $\F\mathcal{A}_{d}$.
In particular, if $n\geq 5$ and $d<2\left\lfloor\frac{n-1}{4}\right\rfloor-1$, then ${\rm SL}_n(\mathbb{Z})$ has pro- perty~$\F\mathcal{A}_{d}$.
\end{corollary}

\begin{remark}
$ $
\begin{enumerate}
\item[$(i)$] \textsc{Bridson} proved in personal communication a slightly better bound for \Aut\ for property $\F\mathcal{A}_{d}$, namely the bound $\lfloor\frac{2n}{3}\rfloor$. 
\item[$(ii)$] There exists an upper bound on the dimension $d$ such that \Aut\ can have property  $\F\mathcal{A}_{d}$.
Consider the symmetric space ${\rm P}_n(\mathbb{R})$ of positive definite real~$n\times n$ matrices. This space is a complete \CAT\ space of dimension $\frac{1}{2}n(n+1)$. The group ${\rm GL}_n(\R)$ acts by isometries on this space via $X\mapsto AXA^t$, where $A\in{\rm GL}_n(\mathbb{R})$, $X\in{\rm P}_n(\mathbb{R})$ and $t$ denotes the transposition. Therefore we have
\[
\Autn\rightarrow{\rm GL}_n(\mathbb{R})\rightarrow{\rm Isom}({\rm P}_n(\mathbb{R}))
\]
and this action is fixed point free. 
\end{enumerate}
\end{remark}

\section*{Acknowledgments} 
The author would like to thank L.~Kramer for his support and for his helpful comments on this article. Many thanks to Matthias Blank for his helpful suggestions and to Daniel Skodlerack for his comments on Theorem B. Finally, the author thanks the referee for many useful suggestions.

This work will be a part of the author's PhD thesis. During the preparation
of this paper the author was supported by the Studienstiftung des deutschen Volkes, Telekom Stiftung and the SFB 878 in M\"unster.

\end{document}